\documentclass[a4paper,12pt,epsfig]{article}
\usepackage{amsmath,amssymb,amsthm,amscd}
\usepackage{graphics}

\usepackage{url}
\newtheorem*{thm}{Theorem}

\theoremstyle{definition}

\theoremstyle{remark}

\newcommand{\s}{\mathrm{s}}

\newcommand{\dist}{\mathop{\rm dist}\nolimits}

\newcommand{\cg}{\mathop{\rm CG}\nolimits}

\title {Towards a proof of the 24--cell conjecture}

\author {Oleg R. Musin\thanks{This research is partially supported by the NSF grant DMS-1400876 and the RFBR grant 15-01-99563.}}

\begin{document}
\date{}
\maketitle

\begin{abstract} This review paper is devoted to the problems of sphere packings in 4 dimensions. The main goal is
to find reasonable approaches for solutions to problems related to densest sphere packings in 4-dimensional
Euclidean space. We consider two long--standing open problems:
the uniqueness of maximum kissing arrangements in 4 dimensions and the 24--cell conjecture. Note that a
proof of the 24--cell conjecture also proves that the lattice packing $D_4$ is the densest sphere packing in 4 dimensions.
\end{abstract}

\medskip


\medskip

\noindent{\bf Acknowledgment.} I wish to thank  W\"oden Kusner  for helpful discussions, comments and especially for detailed references on the Fejes T\'oth/Hales method  (Subsections 2.1 and 2.2).

\section{Introduction}

This paper devoted to the classical problems related to sphere packings in four dimensions.

The sphere packing problem asks for the densest packing of ${\mathbb R}^n$ with unit balls. 
Currently, this problem is solved only for dimensions $n=2$ (Thue [1892, 1910] and  Fejes T{\'o}th [1940], see \cite{kbor04,cohn17b,fej53,zong} for detailed accounts and bibliography), $n=3$ (Hales and Ferguson [42--52]), $n=8$ (Viazovska \cite{viaz17}) and $n=24$ (Cohn et. al. \cite{cohn24}).  

% \cite{hal92,hal97a,hal97b,hal97d,hal00,hal05,hal06a}

In four dimensions, the old conjecture states that a sphere packing is densest when spheres are centered at the points of lattice $D_4$, i.e.  the highest density  $\Delta_4$  is $\pi^2/16$, or equivalently the highest  center density is $\delta_4=\Delta_4/B_4=1/8$.  
For lattice packings, this conjecture was proved by Korkin and Zolatarev in 1872 \cite{KZ1872, KZ1877}. 
Currently, for general sphere packings the best known upper bound for $\delta_4$ is $0.130587$ \cite{LOV14}, a slight improvement on the Cohn--Elkies bound  of $\delta_4<0.13126$ \cite{coh03}, but still nowhere near sharp. 

Consider the Voronoi decomposition of any given packing $P$ of unit spheres in ${\mathbb R}^4$.
\smallskip 

\noindent {\bf The 24--cell conjecture.} {\it The minimal volume of any cell in the resulting Voronoi decomposition of $P$ is at least as large as the volume of a regular 24--cell circumscribed to a unit sphere.}

\smallskip
\noindent Note that a proof of the 24--cell conjecture also proves that $D_4$ is the densest sphere packing in 4 dimensions.

The maximum possible number of non-overlapping unit spheres that can touch a unit sphere in $n$ dimensions is called the kissing number. 
The problem for finding kissing numbers $k(n)$ is closely connected to the more general problems of finding bounds for spherical codes and sphere packings \cite{boy12}. 
Currently, only six kissing numbers are known: $k(1)=2$, $k(2)=6$ (these two are trivial), $k(3)=12$ (some incomplete proofs appeared in 19th century and Sch\"utte and van der Waerden
\cite{sch53} first gave a detailed proof in 1953, see also \cite{lee56,mae01,mae07}), $k(4)=24$ (finally proved in 2003, see  \cite{Mus1} and \cite{mus08a}), $k(8)=240$ and $k(24)=196560$ (found independently in 1979 by Levenshtein \cite{lev79} and Odlyzko--Sloane \cite{odl79}).   
Moreover, Bannai and Sloane \cite{BS} %(= \cite[Chapter 14]{CS}) 
proved that the maximal kissing arrangements in dimensions 8 and 24 are unique up to isometry.  
In dimension 4 the uniqueness of the maximal kissing arrangement is conjectured but not yet proven. 

The main goal is to find reasonable approaches for solutions to problems related to densest sphere packings in 4-dimensional Euclidean space. 
As a basis for this research, we will consider two long-standing open problems: 
the uniqueness of maximum kissing arrangements in 4 dimensions and the 24-cell conjecture.

The paper also considers the following related problems in 4 dimensions: 
the enumeration of optimal and critical spherical configurations of $N$ points for small $N$, which subsumes the study of optimal spherical codes and packings in $\mathbb{S}^3$;  the enumeration of all spherical and Euclidean 2-distance sets;  and the duality gap for LP and SDP bounds.

The initial aims of the project are to examine this duality gap in global LP and SDP problems while simultaneously analyzing the combinatorial structures coming from candidate counterexamples to the 24-cell conjecture defined by unweighted Voronoi cells, as well as those coming from augmented density functionals. 
The goal is to reduce the global packing problem to a local problem in the spirit of Fejes T\'oth and to use a combination of counterexample elimination and SDP techniques to make the local computation tractable.  
The 24-cell conjecture is the most direct reduction and the driving force behind this project.

Our ideas for research problems and preliminary findings are presented in subsequent sections of the paper.

\section{Overview of methods for sphere packing problems} 

\subsection{Problems and methods}

To introduce our problems, let $\mathcal{C} = \{x_1,\dots,x_M\}\subset \mathbb{S}^{d-1}$ be a subset of points on the sphere in $\mathbb{R}^d.$ 
We will call $\mathcal{C}$ a spherical $\varphi$-code if the angular distance between any two points of $\mathcal{C}$ is not less than $\varphi.$ 
By $A(d,\varphi)$ we denote the maximum cardinality of a $\varphi$-code in $\mathbb{S}^{d-1}$. 
For $\varphi=\pi/3$ the problem of finding $A(d,\pi/3)$ is the kissing number problem (see extensive literature). 
For $d=3$, the problem of finding the maximal $\varphi$ such that $A(3,\varphi)\ge n$ for given $n$ is the Tammes problem.

There are three classic methods used for finding the densest sphere packings in metric spaces. 
The local density method goes back to Fejes T{\'o}th \cite{fej43b}, who calculated the maximal density of a sphere packing by considering a triangle with vertices at circle centers and calculating the maximal part of the triangle occupied by circles.
%?????????? 
Coxeter \cite{cox63a} applied this approach to spheres in higher dimensions and conjectured the general upper bounds on $A(d,\varphi)$ by calculating the volume of a regular simplex with edges of angular length $\varphi$ and spatial angle measures at its vertices. 
B{\"o}r{\"o}czky \cite{bor78} verified Coxeter?s conjecture for spaces of constant curvature and thereby proved the Coxeter bound. 
The Coxeter bound is tight for the $600$--cell and therefore $A(4,\pi/5)=120$.

%????????????
Similar ideas are often applied to sphere packing problems of $\mathbb{R}^d$ and particularly to the famous Kepler conjecture \cite{kep11}. 
The choice of partition used for the local approach is especially important in this case. 
Fejes T{\'o}th \cite{fej53} and Hsiang in his unconvincing approach \cite{hsi93a} suggested to use averaging of Voronoi cell densities. 
Hales proposed a local density inequality based on Delaunay triangulations \cite{hal92}, then he formulated inequalities on a ``hybrid'' between Delaunay and Voronoi cells \cite{hal06d}. Finally, the local density inequality of Hales and Ferguson giving the solution to the Kepler conjecture uses the triangulation of space into non-Delaunay triangles \cite{hal06a}.
A simplified method for the formal proof of Kepler uses a hybridization and truncation method introduced by Marchal \cite{mar}.

Also of importance are structural results where local constraints force a global behavior.  In packings of the plane, it is a straight forward observation that the condition that the contact graph is six regular forces a lattice structure; in three dimensions, 12-regularity forces a Barlow packing \cite{hal12}.  It turns out that the determination of the exact value of the kissing radius for 13 points in \cite{mus12a} allows for an alternative proof of this structure theorem in three dimensions \cite{bor15}.  In particular, the exact values found in \cite{mus12a,MTT14} are also tight enough to pass through a series of geometric inequalities and constrain the discrete structure of the spherical Delaunay polytope enough to determine that it must be a rhombic dodecahedron or a triangular orthobicupola.  The same methods of spherical geometry are applicable in the 4-dimensional problem and could be used to reduce the complexity of the case analysis for all the problems we wish to address.

%Irreducible graphs.

%LP and SDP approach.

%\subsection{(a) Densest sphere packings and the kissing number problem}
%\subsection{(b) Methods}
%\subsection{(c) On the uniqueness of kissing arrangement in four dimensions}
%\subsection{(d) The 24-cell conjecture}

%Consider the Voronoi decomposition of any given packing $P$ of unit spheres in ${\mathbb R}^4$.

%\medskip

%\noindent {\bf Conjecture.} {\it The minimal volume of any cell in the resulting Voronoi decomposition of $P$ was at least as large as the volume of a regular 24--cell circumscribed to a unit sphere.}

%\medskip

%The 24--cell conjecture yields the densest packing conjecture in four dimension.

%\medskip

%How do we propose to attack the 24-cell conjecture? 
%
%\noindent I. {\it  By enumeration of neighboring irreducible contact graphs.} 
%\noindent II. {\it By positive definite functions in ${\mathbb R}^4$.} 
%\noindent III. {\it By combination I and II} 

%\subsection{Graphs of sphere packings}

\subsection{Fejes T\'oth/Hales Method: Kepler and Dodecahedral conjectures}
The solution to \emph{the Kepler conjecture}, completed by Hales and Ferguson in \cite{hal06a}, roughly followed an outline proposed by L\'aszl\'o Fejes T\'oth in \cite{fej53}.  
In the same book, Fejes T\'oth linked the kissing problem of Newton and Gregory and the problem of minimal volume configurations in \emph{the Dodecahedral conjecture}, now also a theorem of Hales and McLaughlin \cite{hal10a}.

\begin{thm}[Kepler conjecture: Hales and Ferguson]
There is no packing of $\mathbb{R}^3$ by congruent balls with a density greater than that achieved by the $A_3$ lattice.
\end{thm}

\begin{thm}[Dodecahedral conjecture: Hales and McLaughlin]
There is no Voronoi cell in any unit sphere packing with a volume that less than %that  defined by 12 balls kissing a central ball at the centers of the faces of a regular dodecahedron.
the volume of a regular dodecahedron circumscribed to a unit sphere.
\end{thm}

The final proof method of the Kepler conjecture differs from the strategy proposed by L. Fejes Toth in several ways; it even had to be adapted and evolve over the course of the solution -- but the philosophy is one that we will outline here.  
It is in many ways similar to proposed attacks on the kissing problem and the problem of best lattice packing. Many of such problems are known to be solvable algorithmically.  
The lattice case is solvable via an algorithm due to Voronoi \cite{schur09}, and many geometric problems may be subsumed into the much broader class of optimization over semi-algebraic sets; as long as there are algebraic constraints and objective, the Tarski-Seidenberg algorithm applies.
However, such an approach is intractable in all but the simplest of cases.  
This brute force computational method is not how these problems are generally approached, even if the apparent size of the case analysis makes it appear this way -- the case analyses for the Kepler and Dodecahedral conjectures are massive reductions of the semi-algebraic problem.

In these settings, it is possible to attach relatively simple combinatorial structures to configurations and, via a much smaller enumeration and classification, arrive at a solution.  
The case of arbitrary packings is not generally known to be a finite problem.  
By periodic approximation, it is know that a counterexample to the Kepler conjecture would force the existence of a \emph{finite counterexample}, however, the proof depends on disproving the existence of such. 
To succeed, there must be an auxiliary function that can be attached to the density functional that eliminates all configurations that are sufficiently large. 
In practice, the cutoff is fairly small; there is a decomposition and auxiliary function that works for packings with centers constrained to be with a ball of radius 2.52 (relative to a packing of balls with unit radius).  
The final proof by Hales and Ferguson formed one of the longest proofs in mathematics (a simplification \cite{hal12b} was required to outline the formal verification project, since completed \cite{hal15}).
It was only known \emph{a posteriori} that the Kepler problem was solvable by considering bounded clusters; it was not a given that such an analysis would be guaranteed to terminate.
For the Dodecahedral conjecture, there is a similar cutoff; since it is a local problem, such a cutoff clearly exists, but there is a tradeoff due to the constraint on the decomposition of space -- the auxiliary function must deeply respect the volume of the Voronoi cell.  
The final proof by Hales and McLaughlin initially depended heavily on the machinery developed in the proof of the Kepler conjecture; many of the constructions are directly transferable.  Constructions that are not turn out at least to be transferable by ``close analogy''.

When considering the density of arrangements of spheres in $\mathbb{R}^3$, there are counterexamples to the local optimality of the $A_3$ lattice for obvious decompositions of space.  For example, the Dodecahedral conjecture arose from this fact: 
The Voronoi domain (a regular dodecahedron) of a central sphere, kissing 12 others at the Tammes optimizer is also the volume minimizer among Voronoi cells in sphere packings. 
This is \emph{not} the configuration found in the optimal $A_3$ lattice.  

This might be considered the initial observation in the proof: there might be counterexamples to the global solution.
As observed above, by choosing a clever decomposition of space and attaching an auxiliary function that defines a method of borrowing volume, the existence of a counterexample becomes a finite problem.  
To all configurations, and in particular, to all such potential counterexamples, a combinatorial object can be attached, dependent on the decomposition and auxiliary functions.  
To be attached to a counterexample, these combinatorial objects must satisfy some topological or combinatorial properties: these are the \emph{tame} graphs or hypermaps that must be enumerated.  
In the case of the pure Voronoi decomposition, such a tame graph exists attached to a realizable geometric configuration (the Dodecahedral configuration) which blocks the proof by that particular method. 
 But there exists another decomposition (in fact it appears that there is a large family) when no tame graphs can be realized except the cells associated to the conjectured best packings.  This eliminates all candidate counterexamples and proves the conjecture.
The strategy for the proof of the Dodecahedral conjecture proceeds similarly, but not identically.  In particular, the two proofs illustrate that \emph{tameness} is dependent on the problem; the set of tame graphs must classify the problem and also be constrained enough to enumerate.  In fact, many classical proofs and bounds for packing problems can be placed into this framework

\subsection{Tammes' problem, Fejes T\'oth's method and irreducible contact graphs}

If $N$ unit spheres kiss the unit sphere in ${\Bbb R}^n$, then the set of kissing points is an arrangement on the central sphere such that the (Euclidean) distance between any two points is at least 1. 
This allows us to state the kissing number problem in another way: 
How many points can be placed on the surface of ${\Bbb S}^{n-1}$ so that the angular separation between any two points will be at least $60^{\circ}$?

This leads to an important generalization: 
a finite subset $X$ of ${\Bbb S}^{n-1}$ is called a {\it spherical $\psi$-code} if for every pair  $(x,y)$ of $X$ with $x\ne y$, its angular distance $\dist(x,y)$ is at least $\psi$.

Let $X$ be a finite subset in a metric space $M$.  
Denote 

$$
\psi(X):=\min\limits_{x,y\in X}{\{\dist(x,y)\}}, \mbox{ where } x\ne y.
$$

Let $M={\mathbb S}^2$. 
Denote by $d_N$ the largest angular separation $\psi(X)$ with $|X|=N$ that can be attained in ${\Bbb S}^{2}$, i.e.

$$ 
d_N:=\max\limits_{X\subset{\Bbb S}^2}{\{\psi(X)\}}, \, \mbox{ where } \;  |X|=N.
$$
Consider configurations in ${\mathbb S}^2$ with $\psi(X)=d_N$. 
In other words, {\it how are $N$ congruent, non-overlapping circles distributed on the sphere when the common radius of the circles is as large as possible?}

This question is also known as the problem of the ``inimical dictators'': {\it Where should $N$ dictators build their palaces on a planet so as to be as far away from each other as possible?} The problem was first asked by the Dutch botanist Tammes (see \cite[Section 1.6: Problem 6]{bra05}, who was led to this problem by examining the distribution of pores on the pollen grains of different flowers.

The Tammes problem is presently solved for only a few values of $N$: for $N=3,4,6,12$ by L. Fejes T\'oth \cite{fej43b}; for $N=5,7,8,9$ by Sch\"utte and van der Waerden \cite{sch51}; for $N=10,11$ by Danzer \cite{dan86}; and for $N=24$ by Robinson \cite{rob61}. 
Recently, the problems for $N=13$ and $N=14$ were solved with computer assistance \cite{mus12a, MTT14}. 

The local density method goes back to Fejes T{\'o}th \cite{fej43b}, who calculated the maximal density of a sphere packing by considering a triangle with vertices at circle centers and calculating the maximal part of the triangle occupied by circles. 
He found the following bound:

$$
A(3,\varphi)\le \frac{2\pi}{\Delta(\varphi)}+2,
$$
where 
$$
\Delta(\varphi)=3\arccos{\left(\frac{\cos{\varphi}}{1+\cos{\varphi}}\right)}-\pi,
$$
i.e. $\Delta(\varphi)$ is the area of a spherical regular triangle with side length $\varphi$.

This bound is tight for the Tammes problem for $N=3, 4, 6, 12$, where the configurations are regular triangulations of the sphere. It is also tight asymptotically, since the densest planar circle packing is formed by the regular triangle lattice. For all other cases of the Tammes problem, the Fejes T{\'o}th upper bound can not be tight. 
Robinson \cite{rob61} extended Fejes T\'oth's method and gave a bound valid for all $N$ that is also sharp for $N=24$. 

%Coxeter \cite{cox63a} applied this approach to spheres in higher dimensions and conjectured the general upper bounds on $A(d,\varphi)$ by calculating the volume of a regular simplex with edge angular length $\varphi$ and spatial angle measures at its vertices. B{\"o}r{\"o}czky \cite{bor78} showed that this approach was correct for spaces of constant curvature and thereby proved the Coxeter bound. The Coxeter bound is tight for $120$-cell and therefore $A(4,\pi/5)=120$.

The solutions of all other known cases are based on the investigation of the so-called contact graphs associated with a finite set of points.
For a finite set $X$  in ${\mathbb S}^2$, the {\it contact graph} $\cg(X)$ is the graph with vertices in $X$ and edges $(x,y), \, x,y\in X$, such that $\dist(x,y)=\psi(X)$.  
If the configuration of spherical caps  in ${\mathbb S}^2$ centered in $X$ of diameter $\psi(X)$ is locally rigid, then the graph $\cg(X)$ is said to be irreducible. 
Thus, the study of rigid packings reduces locally to the study of irreducible graphs.

The concept of  irreducible contact graphs was first used by Sch\"utte and van der Waerden to address Tammes' problem \cite{sch51}. 
They used the method also for the solution to the thirteen spheres (kissing number) problem \cite{sch53}. 
In Chapter VI of the  Fejes T\'oth book \cite{fej53},  irreducible contact graphs are considered in greater detail. 
Moreover, in this chapter, solutions for Tammes' problem are conjectured for $N\le 16$, $N=24$ and $N=32$.
The method of irreducible spherical contact graphs was used also \cite{bor83,bor03b,bor03a, bor03,dan86,mus12a, MTT13,MTT14,MT15a} to obtain bounds for the kissing number and Tammes problems. 

The computer-assisted solution of Tammes' problem for $N=13$ and $N=14$ consists of three parts: 
(i) creating the list $L_N$ of all planar graphs with $N$ vertices that satisfy the conditions of \cite[Proposition 3.1]{MTT14}; 
(ii) using linear approximations and linear programming to remove from the list $L_N$ all graphs that do not satisfy the known geometric properties of the maximal contact graphs \cite[Proposition 3.2]{MTT14}; 
(iii) proving that among the remaining graphs in $L_N$ only one is maximal. 

In \cite{mus12b} we considered packings of congruent circles on a square flat torus, i.e., periodic (w.r.t. a square lattice) planar circle packings, with the maximal circle radius. 
This problem is interesting due to a practical reason --- the problem of ``super resolution of images.'' 
We have found optimal arrangements for $N=6$, $7$ and $8$ circles. 
Surprisingly, for the case $N=7$ there are three different optimal arrangements. 
Our proof is based on a computer enumeration of toroidal irreducible contact graphs.

\subsection{LP and SDP methods for sphere packings}

Let $M$ be a metric space with a distance function $\tau.$ 
A real continuous function $f(t)$  %($0\le t \le $ diameter of $M$)
is said to be positive definite (p.d.) in $M$ if
for arbitrary points $p_1,\ldots,p_r$ in $M$, real variables $x_1,\ldots,x_r$, and arbitrary $r$  we have
$$
\sum\limits_{i,j=1}^r {f(t_{ij})\,x_ix_j}\ge 0, \quad t_{ij}=\tau(p_i,p_j),
$$
or equivalently,  the matrix $\bigl(f(t_{ij})\bigr)\succeq0$, where  the sign $\succeq 0$ stands for: ``is  positive semidefinite".

%Let ${\Bbb S}\sp {n-1}$ denote the unit sphere in ${\Bbb R}\sp n$, and let $\varphi_{ij}$ denote the angular distance between points $p_i, p_j.$  

Schoenberg \cite{Scho} proved that: {\em $f(\cos\varphi)$ is p.d. in ${\Bbb S}\sp {n-1}$ if and only if
$\; f(t)=\sum_{k=0}^\infty{f_kG_k^{(n)}(t)}$ with all $f_k\ge 0$.}
 Here $G_k^{(n)}(t)$ 
are the Gegenbauer  polynomials.

Schoenberg's theorem has been generalized by Bochner \cite{boc41} to more general spaces. Namely, the following fact holds: {\em $f$ is p.d. in a 2-point--homogenous space $M$ if and only if $f(t)$ is a nonnegative linear combination   of the zonal spherical functions  $\Phi_k(t)$} (see details in \cite{del77,Kab}, \cite[Ch. 9]{CS}).

The Bochner - Schoenberg theorem plays a crucial role in Delsarte's linear programming (LP) method for finding bounds for the density of sphere packings on spheres and Euclidean spaces. One of the most exciting applications of Delsarte's method is a solution of the kissing number problem in dimensions 8 and 24. However, 8 and 24 are the only dimensions in which this method gives a precise result. For other dimensions (for instance, 3 and 4) the upper bounds exceed the lower. We have found an extension of the Delsarte method \cite{mus06,mus06b,mus08a} that allows to solve the kissing number problem (as well as the one-sided kissing number problem) in dimensions 3 and 4. This method  is widely used in coding theory and discrete geometry for finding bounds for error-correcting codes,  spherical codes, sphere packings and other packing problems in 2-point--homogeneous spaces ([5--7,26--33,55,62--64,88,94,95] and many others). 
%(\cite{ban09a,BS,bar11,cohn17b,cohn07a,CS, Kab, lev79,lev92,lev98,mus11, odl79,pfe04} and many others).

% (\cite{ban09a,BS,bar11,cohn17b,cohn07a,CS, Kab, lev79,lev92,lev98,mus11, odl79,pfe04} and many others).

Cohn and Elkies developed an analogue of the Delsarte LP bounds \cite{coh03} for sphere packing in ${\mathbb R}^n$. (Note that Gorbachev \cite{gor00} independently obtained similar results.)  Using this method, Cohn and Kumar \cite{cohn09} proved the optimality and uniqueness of the Leech lattice among lattices in dimension 24  (see also \cite{pfe04} for a beautiful exposition).  
 Recent solutions of the densest packing problem in dimensions 8 (Viazovska \cite{viaz17}) and 24 (Cohn et. al. \cite{cohn24}) also rely on the Cohn--Elkies method (see also \cite{cohn17}). 
 
{\em Semidefinite programming (SDP)} is a subfield of convex optimization concerned with the optimization of a linear objective function  over the intersection of the cone of positive semidefinite matrices with an affine space.  Schrijver \cite{sch05}  improved some upper bounds on binary codes using  SDP.
Schrijver's method has been adapted for spherical codes  by Bachoc and Vallentin \cite{bac08a}. Now there are many applications of SDP bounds 
 to spherical codes and sphere packings  (see \cite{bac10,bac09a,bac09b,cohn17b,gla13a,mus08b,mus10a,mus14} and many others).

\section{Kissing arrangement uniqueness and the 24-cell conjecture}

\subsection{Contact graphs in four dimensions}
Much of the machinery needed for the analysis of kissing configurations in $\mathbb{R}^4$ can be developed from off-the-shelf components that have been extensively studied and documented; 
such kissing configurations reduce to configurations of spheres in $\mathbb{S}^3$ and the analogous problem in $\mathbb{R}^3$ is the elimination of the local counterexamples to the Kepler problem, as discussed previously.  
The combinatorial methods need to be modified to a spherical geometry, but the considerations would be simpler than those required to determine the best packing in $\mathbb{R}^3$, as there are no candidate configurations that are not in contact with a central sphere -- 
it suffices to start by considering exactly the nearest neighbor contact graphs (or complexes) and characterize them as irreducible with geometric conditions on points and edges (and also higher dimensional cells).  A generally agreed upon condition for irreducibility is local rigidity, there is no shift of a single geometric vertex can increase edge lengths of the contact graph \cite{sch51}. 
%In the case of the 24-cell,  the facets are regular octahedra, 
 Further geometric constraints on the degree and diameter of candidate irreducible or ``tame'' graphs in this context make it a reasonable first approach to the uniqueness conjecture and provide a method to address other kissing configuration problems in $\mathbb{R}^4$.  In particular, we think that more constraints on the combinatorial structure for some optimal packing configurations that considered in recent papers \cite{BBGA17,BG17,BHM17,kus16} can be useful for a proof of  the uniqueness of the 24--cell sphere configuration  in $\mathbb{R}^4$.   %Then, for a sufficiently constrained problem, the uniqueness the solution would follow from the stability of solution to 

\subsection{SDP and uniqueness of the kissing arrangement} Odlyzko and  Sloan \cite{odl79} show that the LP upper  bound for the kissing number $k(4)$ is 25.558... We proved that  $k(4)<24.865$ \cite{mus08a}. 

Denote by $s_d(n)$ the optimal SDP bound on $k(n)$ of degree $d$ \cite{mv10}.   
In the following table it is shown that this minimization problem is a semidefinite program and that every upper bound on $s_d(4)$ provides an upper bound for the kissing number in dimension 4. 

\begin{itemize}

\item $s_7(4) <24.5797$  -- Bachoc \&  Vallentin \cite{bac08a};

\item $s_{11}(4) < 24.10550859$ --   Mittelmann \& Vallentin \cite{mv10};
 
\item $s_{12}(4)< 24.09098111$ \cite{mv10};

\item $s_{13}(4)< 24.07519774$ \cite{mv10};

\item $s_{14}(4) <24.06628391$ \cite{mv10};

\item $s_{15}(4) <24.062758$ -- Machado \&  de Oliveira Filho \cite{mac16};
 
\item $s_{16}(4) <24.056903$  \cite{mac16}.

\end{itemize}

Clearly, the numbers $s_d(n)$ form a monotonic decreasing sequence in $d$. Perhaps this sequence for $n=4$ approaches 24.  If there is a $d$ such that $s_d(4)=24$, then we think it will be possible to prove the uniqueness theorem by a similar way as for dimensions 8 and 24. 

However, since $s_d(4)$ is close to 24 the correspondent polynomial $f$ gives some inequalities for the distances distribution  (see \cite[Theorem 5.4]{mus14}). Moreover, that yields certain constraints for the contact graphs. Therefore, it can help to reduce the list of possible irreducible contact graphs of the kissing arrangements in four dimensions.  

Another interesting possibility is to find an SDP version of our Theorem 1 in \cite{mus08a}. For $n=4$ and $t_0=0.6058$ we have $k(4)\le \max\{h_m\},$ $1\le m\le6$ \cite[Corollary 3]{mus08a}. Then using SDP method we certainly have less $t_0$ and therefore less number of possible configurations. It can also lead to a proof of the uniqueness theorem.

\subsection{The 24--cell conjecture}  The 24-cell conjecture in particularly states that, in the optimal case, all $N$ neighboring spheres touch the central sphere. Since $k(4)=24$, it can be only for $N\le 24$. In order to eliminate the case $N\ge25$ and consider the case $N\le24$, we can use the ideas from the two previous subsections. For this, we can generalize the SDP method for points in ${\mathbb R}^n$ and to apply certain isoperimetric inequalities for polyhedrons.  

In \cite[Sec. 4]{mus14} we define positive--definite (p.d.) functions in ${\mathbb R}\sp n$ 
$H_k^{(n,m)}(t,x,y,{\bf u},{\bf v})$, where $0\le m\le n-2,$ $t,x,y\in {\Bbb R},$ ${\bf u},{\bf v}\in {\Bbb R}\sp m$. Note that for $x=y=1$, $H_k^{(n,0)}$ is the Gegenbauer polynomial $G_k^{(n)}(t)$, and for $m=1$ that is the multivariate Gegenbauer  polynomial   $S_k^n(t,u,v)$ first defined by Bachoc and Vallentin \cite{bac08a}. 

Let $H_k:=H_k^{(4,1)}$. Then $H_k$ is a p.d. polynomial in five variables $t, x, y, u, v$.   If $p_1,\ldots, p_N$  in ${\Bbb R}\sp 4$ are centers of unit neighboring spheres with the central sphere centered at the origin, then $t=\langle p_i,p_j\rangle,$ $x=|p_i|^2$ and $y=|p_j|^2$, see \cite[Theorem 4.1]{mus14}.  Since all $|p_i|$ are close to 1, all $H_k$ are close to $S_k^4$. Thus, perhaps we can have similar bounds as for the spherical case and in particular $N\le24$.   

%Let $m=1$ and $n=4$. Then $H_k:=H_k^{(n,m)}$ is a polynomial of degree $k$ in five variables $t,x,y,u,v$. 

%Theorem 4.1 \cite{mus14} states that for any set of points $p_1,\ldots, p_N$  in ${\Bbb R}\sp n$  the matrix $\bigl(g_{ij}\bigr)$ is positive semidefinite, where $$g_{ij}=H_k^{(n,m)}(\langle p_i,p_j\rangle,|p_i|^2,|p_j|^2,p_i,p_j).$$

\subsection{Dimension reduction}
 Here we consider more difficult problems for which we do not have a ready approach in mind, but which we still wish to analyze. 
 If we increase the degree of the polynomial, then the dimension of the SDP problem rapidly increases.
  It seems to us that it is possible to apply the methods of combinatorial topology, namely fixed-point theorems for reducing the dimension of the corresponding SDP problems. 

One of the successful implementations of this approach is a paper by  Bondarenko, Radchenko, and Viazovska \cite{BRV} on spherical $t$--designs.  (A {\em spherical $t$--design} is a finite set of $N$ points on  ${\Bbb S}^d$  such that the average value of any polynomial $f$ of degree $t$ or less on the set equals the average value of $f$ on the whole sphere.)
They proved the conjecture of Korevaar and Meyers:

\medskip
 
\noindent{\em For each $N\ge c_dt^d$ there exists a spherical $t$--design in the sphere ${\Bbb S}^d$ consisting of $N$ points, where $c_d$ is a constant depending only on $d$.}

\medskip

\noindent One of the most important steps in their proof is based on the lemma that follows from the Brouwer fixed--point theorem.

Note that our topological and topological combinatorics papers \cite{mus12c, MusSpT, MusQ, MusS, MusH, MusKKM, MusVo} are particularly motivated by optimal sphere packing problems.
 
\section{Related research problems}

\subsection{Optimal spherical codes in four dimensions}
From the perspective of spherical codes, it is a shame that so little is known about kissing configurations in higher dimensions. 
 Fejes T\'oth-type inequalities imply that the regular simplex is an optimizer in all dimensions, as well as the octahedron and icosahedron in dimension 3, and the 600-cell in dimension 4.  Otherwise, in dimension 4, the only other spherical codes known to be optimal are for configurations with fewer than 8 points and for exactly 10.  
 It is conjecture that the 9-point configuration is of the form \{1,4,4*\}, one point in the pole and eight vertices of two simplexes, twisted relative to reach other, on ``spheres of latitude".  
 Such problems seem approachable by methods of irreducible graphs related to the classifications that would be required to address the 24-cell conjecture and would yield insight into the discrete geometry of $\mathbb{S}^3$. A further generalization of Robinson's tightening of the Fejes T\'oth inequality would also be of great interest.

\subsection{Maximum contact packings in four dimensions}
The {\em spherical kissing number $k_S(d,\theta)$} is the maximum number of disjoint spherical caps of angular diameter $\theta$ in ${\mathbb S}^{d}$ that can be arranged so that all of them touch one spherical cap of the same diameter. Denote by $k_S(d)$ the maximum value of $k_S(d,\theta)$. In fact,  $k_S(d,\theta)=k_S(d)$ for $\theta\to0$ and $k_S(d)=A(d,\varphi)$ if $\varphi<\pi/3$ is close to $\pi/3$. 

Currently, spherical kissing numbers are known only for $d\le3$. Namely, $k_S(1)=2,$ $k_S(2)=5$ and $k_S(3)=12$. Our conjecture is that $$k_S(4)=22.$$

The  {\em contact graph} of an arbitrary finite packing $P$ of unit spheres in ${\mathbb R}^{d}$  is the graph whose vertices correspond to the packing spheres and whose two vertices are connected by an edge if the corresponding two packing spheres touch each other.   
Denote by $c(n,d)$ the maximum number of touching pairs in packings $P$  of cardinality $n$. In other words, $c(n,d)$ is the maximum number of edges of the contact graph of a packing  of $n$ unit spheres in ${\mathbb R}^{d}$. 

It is clear that 
$$
c(n,d)<\frac{1}{2}k(d)\,n. 
$$
In 1974 Harborth \cite{Harb} proved that 
$$
c(n,2)=\lfloor 3n-\sqrt{12n-3}\rfloor. 
$$
There are only particular results for higher dimensions   \cite{bezdek2012,bezdek2013,bezdek2016}. 

 Denote by $\s(n,d)$ the maximum number of touching pairs in packings  of ${\mathbb S}^{d-1}$  by $n$ congruent spherical caps.  We have 
$$
c_S(n,d)\le\frac{1}{2}k_S(d)\,n. 
$$

Robinson and Fejes T\'oth for $d=2$  found all cases when the equality holds (see \cite{musinLFT} for references and other results).   It is interesting to solve this problem in four dimensions. 

%\subsection{Two--distance sets in four dimensions}

 \subsection{Two-distance sets in four dimensions} 
 A set $S$ in Euclidean space ${\Bbb R}^d$ is called a {\it two-distance set}, if there are two distances $a$ and $b$, and the
distances between pairs of points of $S$ are either $a$ or $b$. If a two-distance set $S$ lies in the unit sphere ${\Bbb S}^{d-1}$,
then $S$ is called a {\it spherical two-distance set.} 

Let $G$ be a graph on $n$ vertices. Consider a {\em Euclidean representation of $G$} in ${\Bbb R}^d$ as a two-distance set. In other words, there are two positive real numbers $a$ and $b$ with $b\ge a>0$ and an embedding $f$ of the vertex set of $G$ into  ${\Bbb R}^d$ such that 
$$ \dist(f(u),f(v)):=\left\{
\begin{array}{l}
a \; \mbox{ if } uv \mbox{ is an edge of } G\\
b \;  \mbox{ otherwise}
\end{array} 
\right.
$$

We will call the smallest $d$ such that $G$ is representable in ${\Bbb R}^d$ the {\em Euclidean representation number} of $G$ and denote it by $\dim_2(G)$. Let $G$ be a simple graph on $n$ vertices. It is clear that $\dim_2(G)\le n-1$. Einhorn and Schoenberg \cite{ES}  proved the following theorem: 

\medskip

\centerline{\em $\dim_2(G)= n-1$ if and only if\, $G$ is a disjoint union of cliques.}

\medskip

Denote by $\Sigma_n$ the number of all two-distance sets with $n$ vertices in ${\Bbb R}^{n-2}$. Then Einhorn--Schoenberg's theorem yields 
 $$
\Sigma_n = \Gamma_n - p(n),
$$
where $\Gamma_n$ is the number of all simple undirected graphs and $p(n)$ is the number of unrestricted partitions of $n$.

Einhorn and Schoenberg \cite{ES} enumerated all two-distance sets in dimensions two and three. In other words, they enumerated all graphs $G$ with $\dim_2(G)=2$ and $\dim_2(G)=3$. This problem in dimension four is still open. 

Einhorn--Schoenberg's theorem  gives a complete enumeration of two distance-sets in ${\mathbb R}^4$ of cardinality $n\le6$. In particular, since $\Gamma_6=156$ and $p(6)=11$, we have $\Sigma_6=145.$  

Lison{\v{e}}k \cite{lisonek1997} proved that the maximum cardinality of two-distance sets in ${\mathbb R}^4$ is 10. Moreover, this representation is unique up to similarity. It remains to solve the problem for $n=7,$ 8 and 9.

In \cite{mus2dist}, we consider the spherical representation number of $G$.  We give exact formulas for this number using multiplicities of polynomials that are defined by the Caley--Menger determinant. We think that using this method can enumerated all spherical two-distance sets in four dimensions.

 \medskip

 \medskip

 \medskip

O. R. Musin, 
%IITP RAS, Bolshoy Karetny per. 19, Moscow, 127994, Russia\\
%and\\
University of Texas Rio Grande Valley, School of Mathematical and Statistical Sciences,  One West University Boulevard,  Brownsville, TX, 78520.\\
{\it E-mail address:} oleg.musin@utrgv.edu


\medskip

\medskip 

\begin{thebibliography}{100}
{\renewcommand\baselinestretch{1.03}\normalsize

%\bibitem{ako11}
%A.~Akopyan, A.~Glazyrin, O.~R. Musin, and A.~Tarasov.
%\newblock The extremal spheres theorem.
%\newblock {\em Discrete Math.}, 311(2-3):171--177, 2011.

%\bibitem{and99}
%G.~E. Andrews, R.~Askey, and R.~Roy.
%\newblock {\em Special functions}, volume~71 of {\em Encyclopedia of
%  Mathematics and its Applications}.
%\newblock Cambridge University Press, Cambridge, 1999.

%\bibitem{ans04}
%K.~M. Anstreicher.
%\newblock The thirteen spheres: a new proof.
%\newblock {\em Discrete Comput. Geom.}, 31(4):613--625, 2004.

%\bibitem{bac06}
%C.~Bachoc.
%\newblock Linear programming bounds for codes in {G}rassmannian spaces.
%\newblock {\em IEEE Trans. Inform. Theory}, 52(5):2111--2125, 2006.

\bibitem{bac10}
C.~Bachoc, D.~C. Gijswijt, A.~Schrijver, and F.~Vallentin.
\newblock Invariant semidefinite programs.
\newblock In {\em Handbook on semidefinite, conic and polynomial optimization},
  volume 166 of {\em Internat. Ser. Oper. Res. Management Sci.}, pages 219--269. Springer, New York, 2012.

\bibitem{bac08a}
C.~Bachoc and F.~Vallentin.
\newblock New upper bounds for kissing numbers from semidefinite programming.
\newblock {\em J. Amer. Math. Soc.}, 21(3):909--924, 2008.

\bibitem{bac09a}
C.~Bachoc and F.~Vallentin.
\newblock Optimality and uniqueness of the {$(4,10,1/6)$} spherical code.
\newblock {\em J. Combin. Theory Ser. A}, 116(1):195--204, 2009.

\bibitem{bac09b}
C.~Bachoc and F.~Vallentin.
\newblock Semidefinite programming, multivariate orthogonal polynomials, and
  codes in spherical caps.
\newblock {\em European J. Combin.}, 30(3):625--637, 2009.

%\bibitem{ban05}
%E.~Bannai and E.~Bannai.
%\newblock A note on the spherical embeddings of strongly regular graphs.
%\newblock {\em European J. Combin.}, 26(8):1177--1179, 2005.

%\bibitem{ban06}
%E.~Bannai and E.~Bannai.
%\newblock On primitive symmetric association schemes with {$m_1=3$}.
%\newblock {\em Contrib. Discrete Math.}, 1(1):68--79 (electronic), 2006.

\bibitem{ban09a}
E.~Bannai and E.~Bannai.
\newblock A survey on spherical designs and algebraic combinatorics on spheres.
\newblock {\em European J. Combin.}, 30(6):1392--1425, 2009.



\bibitem{BS}
E. Bannai and N. J. A. Sloane, Uniqueness of certain spherical codes, {\em Canadian J. Math.,}  {33}:437--449, 1981.

\bibitem{bar11}
A.~Barg and O.~R. Musin.
\newblock Bounds on sets with few distances.
\newblock {\em J. Combin. Theory Ser. A}, 118(4):1465--1474, 2011.

%\bibitem{barg13}
%A. Barg and W.-H. Yu, New bounds for spherical two-distance set, {\em Experimental Math.,}  22:
 %18--194, 2013.
 
% \bibitem{baryshnikov2014min}
%Y~Baryshnikov, P~Bubenik, and M~Kahle.
%\newblock Min-type {M}orse theory for configuration spaces of hard spheres.
%\newblock {\em International Mathematics Research Notices}, 2014(9):2577--2592,
  %2014.
 
 \bibitem{bezdek2012}
K. Bezdek.
\newblock  Contact numbers for congruent sphere packings in Euclidean 3-space.
\newblock {\em  Discrete Comput. Geom.}, 48(2): 298--309, 2012.

\bibitem{bezdek2013}
K. Bezdek and S. Reid.
\newblock  Contact graphs of unit sphere packings revisited.
\newblock {\em J. Geom.}, 104(1): 57--83, 2013.

\bibitem{bezdek2016}
K. Bezdek and M. A. Khan. 
\newblock  Contact numbers for sphere packings.
\newblock arXiv:1601.00145, 2016.

\bibitem{boc41}
S.~Bochner.
\newblock Hilbert distances and positive definite functions.
\newblock {\em Ann. of Math. (2)}, 42:647--656, 1941.


\bibitem{BRV}
A. Bondarenko, D. Radchenko  and M. Viazovska. 
\newblock On optimal asymptotic bounds for spherical designs.
 \newblock {\it Ann. of Math.}, vol. 178,  443--452, 2013. 

\bibitem{bor78}
K.~B{\"o}r{\"o}czky.
\newblock Packing of spheres in spaces of constant curvature.
\newblock {\em Acta Math. Acad. Sci. Hungar.}, 32(3-4):243--261, 1978.

\bibitem{bor83}
K.~B{\"o}r{\"o}czky.
\newblock The problem of {T}ammes for {$n=11$}.
\newblock {\em Studia Sci. Math. Hungar.}, 18(2-4):165--171, 1983.

\bibitem{bor03b}
K.~B{\"o}r{\"o}czky.
\newblock The {N}ewton-{G}regory problem revisited.
\newblock In {\em Discrete geometry}, volume 253 of {\em Monogr. Textbooks Pure
  Appl. Math.}, pages 103--110. Dekker, New York, 2003.
  
 \bibitem{BBGA17}
K.~B{\"o}r{\"o}czky, K.~J. B{\"o}r{\"o}czky, A. Glazyrin, \'A. Kov\'acs.  Stability of the simplex bound for packings by equal spherical caps determined by simplicial regular polytopes, preprint, arXiv:1711.00211 

 \bibitem{BG17}
K.~J. B{\"o}r{\"o}czky, A. Glazyrin, Stability of optimal spherical codes,  arXiv:1711.06012
  
 \bibitem{BHM17}
K.~B{\"o}r{\"o}czky, A. Heppes, E. Makai Jr.,   Densest packings of translates of strings and layers of balls, arXiv:1706.05282

\bibitem{bor03a}
K.~B{\"o}r{\"o}czky and L.~Szab{\'o}.
\newblock Arrangements of 13 points on a sphere.
\newblock In {\em Discrete geometry}, volume 253 of {\em Monogr. Textbooks Pure
  Appl. Math.}, pages 111--184. Dekker, New York, 2003.

\bibitem{bor03}
K.~B{\"o}r{\"o}czky and L.~Szab{\'o}.
\newblock Arrangements of 14, 15, 16 and 17 points on a sphere.
\newblock {\em Studia Sci. Math. Hungar.}, 40(4):407--421, 2003.

\bibitem{bor15}
K.~B{\"o}r{\"o}czky and L.~Szab{\'o}.
\newblock 12--Neighbour packings of unit balls in ${\mathbb E}^3$. 
\newblock {\em Acta Math. Hungar.,} 146 (2): 421--448, 2015.

\bibitem{kbor04}
K.~B{\"o}r{\"o}czky, Jr. Finite Packing and Covering, Cambridge University Press, 2004. 

\bibitem{boy12}
P.~Boyvalenkov, S.~Dodunekov, and O.~Musin.
\newblock A survey on the kissing numbers.
\newblock {\em Serdica Math. J.}, 38(4):507--522, 2012.

\bibitem{bra05}
P.~Brass, W.~Moser, and J.~Pach.
\newblock {\em Research problems in discrete geometry}.
\newblock Springer, New York, 2005.

%\bibitem{bri07}
%G.~Brinkmann and B.~D. McKay.
%\newblock Fast generation of planar graphs.
%\newblock {\em MATCH Commun. Math. Comput. Chem.}, 58(2):323--357, 2007.

%\bibitem{cantarella2006criticality}
%J~Cantarella, J~H~G Fu, R~Kusner, J~M Sullivan, and N~C Wrinkle.
%\newblock Criticality for the {G}ehring link problem.
%\newblock {\em Geometry \& Topology}, 10(4):2055--2115, 2006.

%\bibitem{cas04}
%B.~Casselman.
%\newblock The difficulties of kissing in three dimensions.
%\newblock {\em Notices Amer. Math. Soc.}, 51(8):884--885, 2004.

\bibitem{cohn17}
H.~Cohn. 
\newblock A conceptual breakthrough in sphere packing. 
\newblock {\em Notices Amer. Math. Soc.}, 64:102--115, 2017.

\bibitem{cohn17b}
H.~Cohn. 
\newblock Packing, coding, and ground states, in {\em Mathematics and Materials,} 
\newblock {\em IAS/Park City Mathematics Series}, 24: 45--102, 2017.

\bibitem{coh03}
H.~Cohn and N.~Elkies.
\newblock New upper bounds on sphere packings. {I}.
\newblock {\em Ann. of Math.}, 157(2):689--714, 2003.

\bibitem{cohn07a}
H.~Cohn, J. H. Conway, N.~Elkies and A. Kumar.
\newblock The $D_4$ root system is not universally optimal.
\newblock {\em Experiment. Math.}, 16(3):313--320, 2007.

\bibitem{cohn24}
H. Cohn, A. Kumar, S. D. Miller, D. Radchenko, and M. Viazovska. 
\newblock The sphere packing problem in dimension 24. 
\newblock {\em Ann. of Math.}, 185:1017--1033, 2017.


\bibitem{cohn2007universally}
H~Cohn and A~Kumar.
\newblock Universally optimal distribution of points on spheres.
\newblock {\em Journal of the American Mathematical Society}, 20(1):99--148,
  2007.

%\bibitem{coh07b}
%H.~Cohn and A.~Kumar.
%\newblock Universally optimal distribution of points on spheres.
%\newblock {\em J. Amer. Math. Soc.}, 20(1):99--148, 2007.

\bibitem{cohn09}
H.~Cohn and A.~Kumar.
\newblock Optimality and uniqueness of the Leech lattice among lattices.
\newblock {\em Ann. of Math.}, 170(3):1003--1050, 2009.

\bibitem{coh12}
H.~Cohn and Y.~Zhao.
\newblock Sphere packing bounds via spherical codes.
\newblock {\em Duke Mat, J.,} 163 (10):1965--2002, 2014.

%\bibitem{con05}
%R.~Connelly.
%\newblock Generic global rigidity.
%\newblock {\em Discrete Comput. Geom.}, 33(4):549--563, 2005.

%\bibitem{con08}
%R.~Connelly.
%\newblock Rigidity of packings.
%\newblock {\em European J. Combin.}, 29(8):1862--1871, 2008.

\bibitem{CS}
J.~H. Conway and N.~J.~A. Sloane, Sphere Packings, Lattices, and Groups, New York, Springer-Verlag, 1999 (Third Edition).

\bibitem{cox63a}
H.~S.~M. Coxeter.
\newblock An upper bound for the number of equal nonoverlapping spheres that
  can touch another of the same size.
\newblock In {\em Proc. {S}ympos. {P}ure {M}ath., {V}ol. {VII}}, pages 53--71.
  Amer. Math. Soc., Providence, R.I., 1963.

\bibitem{dan86}
L.~Danzer.
\newblock Finite point-sets on {$\mathbb{S}^2$} with minimum distance as large
  as possible.
\newblock {\em Discrete Math.}, 60:3--66, 1986.

%\bibitem{del72}
%P.~Delsarte.
%\newblock Bounds for unrestricted codes, by linear programming.
%\newblock {\em Philips Res. Rep.}, 27:272--289, 1972.

%\bibitem{del73}
%P.~Delsarte.
%\newblock An algebraic approach to the association schemes of coding theory.
%\newblock {\em Philips Res. Rep. Suppl.}, (10):vi+97, 1973.

%\bibitem{del73a}
%P.~Delsarte.
%\newblock Four fundamental parameters of a code and their combinatorial
%  significance.
%\newblock {\em Information and Control}, 23:407--438, 1973.

\bibitem{del77}
P.~Delsarte, J.~M. Goethals, and J.~J. Seidel.
\newblock Spherical codes and designs.
\newblock {\em Geometriae Dedicata}, 6(3):363--388, 1977.

%\bibitem{dol12a}
%N.~P. Dolbilin, O.~R. Musin, and H.~Edelsbrunner.
%\newblock On the optimality of functionals over triangulations of Delaunay sets.
%\newblock {\em Russian Math. Surv.}, 64(4):781--783, 2012.

%\bibitem{dol12b}
%N.~P. Dolbilin, H.~Edelsbrunner, A.~Glazyrin, and O.~R. Musin.
%\newblock Functionals on triangulations of Delaunay sets.
%\newblock {\em Mosc. Math. J.}, 14(3):491--504, 2014.

% \bibitem{edel17}
 %H. Edelsbrunner, A. Glazyrin, O. R. Musin and  A. Nikitenko.  The Voronoi functional is maximized by the Delaunay triangulation in the plane, awaiting publication  in {\em Combinatorica}, 2017.
 
 \bibitem{ES}
S.~J. Einhorn and I.~J. Schoenberg.
\newblock On euclidean sets having only two distances between points. {I, II}.
\newblock {\em Nederl. Akad. Wetensch. Proc. Ser. A 69=Indag. Math.},
  28:479--488, 489--504, 1966.

\bibitem{fej43b}
L.~Fejes~T{\'o}th.
\newblock {\"U}ber eine {A}bschatzung des k{\"u}rzesten {A}bstandes zweier
  {P}unkt eines auf einer {K}ugelfl{\"u}che liegenden {P}unksystems.
\newblock {\em Jber. Deutschen Math. Verein}, 53:65--68, 1943.

\bibitem{fej53}
L.~Fejes~T{\'o}th.
\newblock {\em Lagerungen in der {E}bene, auf der {K}ugel und im {R}aum}.
\newblock Die Grundlehren der Mathematischen Wissenschaften in
  Einzeldarstellungen mit besonderer Ber{\"u}cksichtigung der
  Anwendungsgebiete, Band LXV. Springer-Verlag, Berlin, 1953.


\bibitem{gla13a}
A.~Glazyrin and O.~R. Musin.
\newblock The price of {SDP} relaxations for spherical codes.
\newblock  in  ``Designs, Codes, Graphs and Related Areas'',  RIMS Kokyuroku, vol. 1889:1--6, 2014.

%\bibitem{glaz16}
%A. Glazyrin and W.-H. Yu, 
%\newblock Upper bounds for $s$--distance sets and equiangular lines,
%\newblock{\em arXiv:1611.09479}, 2016.

\bibitem{gor00}
D.~V. Gorbachev.
\newblock An extremal problem for entire functions of exponential spherical
  type, which is connected with the {L}evenshtein bound for the density of a
  packing of {${\Bbb R}^n$} by balls.
\newblock {\em Izv. Tul. Gos. Univ. Ser. Mat. Mekh. Inform.}, 6(1,
  Matematika):71--78, 2000.

\bibitem{hal92}
T.~C. Hales.
\newblock The sphere packing problem.
\newblock {\em J. Comput. Appl. Math.}, 44(1):41--76, 1992.

\bibitem{hal97a}
T.~C. Hales.
\newblock Sphere packings. {I}.
\newblock {\em Discrete Comput. Geom.}, 17(1):1--51, 1997.

\bibitem{hal97b}
T.~C. Hales.
\newblock Sphere packings. {II}.
\newblock {\em Discrete Comput. Geom.}, 18(2):135--149, 1997.

\bibitem{hal06d}
T.~C. Hales.
\newblock Sphere packings. {III}.
\newblock {\em Discrete Comput. Geom.}, 36(1):71--110, 2006.

\bibitem{hal00}
T.~C. Hales.
\newblock Cannonballs and honeycombs.
\newblock{\em Notices of the AMS}, 47(4):440--449, 2000.

\bibitem{hal05}
T.~C. Hales.
\newblock A proof of the {Kepler} conjecture.
\newblock {\em Ann. of Math.}, 162(3):1065--1185, 2005.



\bibitem{hal15}
T.~C. Hales et. al..
\newblock A formal proof of the {Kepler} conjecture.
\newblock {\em arXiv:1501.02155}, 2015.

\bibitem{hal12}
T.~C. Hales.
\newblock A proof of {F}ejes {T}{\'o}th's conjecture on sphere packings with
  kissing number twelve.
\newblock {\em arXiv:1209.6043}, 2012.

\bibitem{hal12b}
T.~C. Hales.
\newblock {\em Dense Sphere Packings: A Blueprint for Formal Proofs}.
\newblock Cambridge University Press, 2012.

%\bibitem{hal16}
%T.~C. Hales and W. Kusner. 
%\newblock Packings of regular pentagons in the
 % plane
  %\newblock {\em arXiv:1602.07220}, 2016.

\bibitem{hal06a}
T.~C. Hales and S.~P. Ferguson.
\newblock A formulation of the {K}epler conjecture.
\newblock {\em Discrete Comput. Geom.}, 36(1):21--69, 2006.

\bibitem{hal10a}
\newblock T.~C. Hales and S. McLaughlin.
\newblock The dodecahedral conjecture.
\newblock  {\em Journal of the AMS}, 23(2):299--344, 2010.

\bibitem{Harb}
H. Harborth, L\"osung zu Problem 664A, {\em Elem. Math.}, 29 (1974), 14--15.

\bibitem{hsi93a}
W.-Y. Hsiang.
\newblock On the sphere packing problem and the proof of {K}epler's conjecture.
\newblock In {\em Differential geometry and topology ({A}lghero, 1992)}, pages
  117--127. World Sci. Publ., River Edge, NJ, 1993.

%\bibitem{hsi93b}
%W.-Y. Hsiang.
%\newblock On the sphere packing problem and the proof of {K}epler's conjecture.
%\newblock {\em Internat. J. Math.}, 4(5):739--831, 1993.

%\bibitem{hsi01}
%W.-Y. Hsiang.
%\newblock {\em Least action principle of crystal formation of dense packing
 % type and {K}epler's conjecture}, volume~3 of {\em Nankai Tracts in
%  Mathematics}.
%\newblock World Scientific Publishing Co. Inc., River Edge, NJ, 2001.
%\newblock With a foreword by S. S. Chern.

\bibitem{Kab}
G.~A. Kabatiansky and V.~I. Levenshtein, Bounds for packings on a sphere and in space,
% Problemy Peredachi informacii {\bf 14}(1), 1978, 3-25; English translation,
Problems of Information Transmission, 14(1):1--17, 1978.

\bibitem{kep11}
J.~Kepler.
\newblock {S}trena {S}eu de {N}iue {S}exangula ({T}he {S}ix-{C}ornered
  {S}nowflake).
\newblock {\em Published by Godfrey Tampach at Frankfort on Main}, 1611.

\bibitem{kus16}
\newblock R. Kusner, W. Kusner, J. Lagarias, and S. Shlosman.
\newblock The twelve spheres problem. 
\newblock { \em arXiv:1611.10297}, 2016.

\bibitem{KZ1872}
A. Korkin and G. Zolotarev. Sur les formes quadratiques positive quaternaires. {\em Math. Ann.,} 5:581--583, 1872.

\bibitem{KZ1877}
A. Korkin and G. Zolotarev.  Sur les formes quadratiques positives. {\em Math. Ann.,} 11:242--292, 1877.

\bibitem{LOV14}
\newblock D. de Laat,  F. M. de Oliveira Filho and F. Vallentin. 
\newblock Upper bounds for packings of spheres of several radii. 
\newblock {\em Forum Math. Sigma 2,} e23, 2014.

\bibitem{lee56}
J.~Leech.
\newblock The problem of the thirteen spheres.
\newblock {\em Math. Gaz.}, 40:22--23, 1956.

\bibitem{lev79}
V.~I. Levenshtein.
\newblock Boundaries for packings in {$n$}-dimensional {E}uclidean space.
\newblock {\em Dokl. Akad. Nauk SSSR}, 245(6):1299--1303, 1979.

\bibitem{lev92}
V.~I. Levenshtein.
\newblock Designs as maximum codes in polynomial metric spaces.
\newblock {\em Acta Appl. Math.}, 29(1-2):1--82, 1992.
\newblock Interactions between algebra and combinatorics.

\bibitem{lev98}
V.~I. Levenshtein.
\newblock Universal bounds for codes and designs.
\newblock In {\em Handbook of coding theory, {V}ol. {I}, {II}}, pages 499--648.
  North-Holland, Amsterdam, 1998.
  
 \bibitem{lisonek1997}
P.~Lison{\v{e}}k.
\newblock New maximal two-distance sets.
\newblock {\em Journal of Combinatorial Theory, Series A}, 77(2):318--338, 1997. 
  
\bibitem{mac16} 
\newblock F. C. Machado and F. M. de Oliveira Filho.
 \newblock Improving the semidefinite programming bound for the kissing number by exploiting polynomial symmetry, 
 {\em arXiv:1609.05167}, 2016.

\bibitem{mae01}
H.~Maehara.
\newblock Isoperimetric theorem for spherical polygons and the problem of 13
  spheres.
\newblock {\em Ryukyu Math. J.}, 14:41--57, 2001.


\bibitem{mae07}
H.~Maehara.
\newblock The problem of thirteen spheres---a proof for undergraduates.
\newblock {\em European J. Combin.}, 28(6):1770--1778, 2007.

\bibitem{mar}
C.~Marchal.
\newblock Study of Kepler's conjecture: the problem of the closest packing.
\newblock {\em Mathematische Zeitschrift}, 267(3--4):737--765, 2011.

\bibitem{mv10}
H. D. Mittelmann and F. Vallentin,
\newblock High-accuracy semidefinite programming bounds for kissing numbers,
\newblock {\em Experimental Math.}, 9:175--179, 2010.

%\bibitem{mus98}
%O.~R. Musin.
%\newblock {Chebyshev systems and zeros of a function on a convex curve.}
%\newblock  {\em Proc. Steklov Inst. of Math.},  221:236--246, 1998.

\bibitem{Mus1} 
O.~R. Musin.
\newblock The problem of twenty-five spheres,
\newblock {\em Russian Math. Surveys} 58:794--795, 2003.


\bibitem{mus06}
O.~R. Musin.
\newblock The kissing problem in three dimensions.
\newblock {\em Discrete Comput. Geom.}, 35(3):375--384, 2006.

\bibitem{mus06b}
O.~R. Musin, The one-sided kissing number in four dimensions,
{\em Periodica Math. Hungar.}, 53:209--225, 2006.


\bibitem{mus08a}
O.~R. Musin.
\newblock The kissing number in four dimensions.
\newblock {\em Ann. of Math.}, 168:1--32, 2008.


\bibitem{mus08b}
O.~R. Musin.
\newblock Bounds for codes by semidefinite programming.
\newblock {\em Proc. Steklov Inst. Math.}, 263:134--149, 2008.


\bibitem{mus09a}
O.~R. Musin.
\newblock Spherical two-distance sets.
\newblock {\em J. Combin. Theory Ser. A}, 116(4):988--995, 2009.

%\bibitem{mus09b}
%O.~R. Musin.
%\newblock An inverse theorem on equivariant genera.
%\newblock {\em Russian Math. Surveys,}, 64:753--755, 2009.

\bibitem{mus10a}
O.~R. Musin.
\newblock Positive definite functions in distance geometry.
\newblock In {\em European Congress of Mathematics}, pages 115--134. Eur. Math. Soc., Z{\"u}rich, 2010.

%\bibitem{mus11a}
%O.~R. Musin.
%\newblock On rigid Hirzebruch genera.
%\newblock {\em Mosc. Math. J.}, 11:139--147, 2011.


\bibitem{mus12c}
O.~R. Musin.
\newblock Borsuk--Ulam type theorems for manifolds.
\newblock {\em Proc. Amer. Math. Soc.}, 140:2551--2560, 2012.

\bibitem{mus14}
{O. R. Musin}, Multivariate positive definite functions on spheres, in: {\it Discrete Geometry and Algebraic Combinatorics, A. Barg and O. Musin, Editors, AMS Series: Contemporary Mathematics,} vol. {625}:77--190, 2014.


\bibitem{MusSpT} 
O. R. Musin, Extensions of Sperner and Tucker's lemma for manifolds. {\it J. of Combin. Theory Ser. A,} 132:172--187, 2015.

\bibitem{MusQ}
O. R. Musin, Sperner type lemma for quadrangulations. {\it Moscow Journal of Combinatorics and Number theory}, 5:26--35, 2015.  

\bibitem{MusS}
O. R. Musin, Generalizations of Tucker--Fan--Shashkin lemmas. {\it Arnold Math. J.},  2(3):299--308, 2016. 

\bibitem{MusH} 
O. R. Musin, Homotopy invariants of covers and KKM type lemmas.  {\it Algebr. Geom. Topol.}, 16:1799--1812, 2016.

%\bibitem{MusC2} 
 %{{O. R. Musin,} Circle actions with two isolated points. {\it  Math. Notes}, {100(4)}:134--136, 2016.}
 
 \bibitem{MusKKM} 
O. R. Musin, KKM type theorems with boundary conditions. {\em J. Fixed Point Theory Appl.}, {19}:2037--2049, 2017.

\bibitem{musinLFT}
O.~R. Musin.
\newblock  Five Essays on the Geometry of L\'aszl\'o Fejes T\'oth. 
\newblock  {\em arXiv:1604.02776}

\bibitem{mus2dist}
O.~R. Musin. Graphs and spherical two-distance sets. {\em arXiv:1608.03392}

\bibitem{mus12b}
O.~R. Musin and A.~V. Nikitenko.
\newblock Optimal packings of congruent circles on a square flat torus.
\newblock  {\em Discrete Comput. Geom.},  {55(1)}:1--20, 2016. 

\bibitem{mus11}
O.~R. Musin and H.~Nozaki.
\newblock Bounds on three- and higher-distance sets.
\newblock {\em European J. Combin.}, 32(8):1182--1190, 2011.

\bibitem{mus12a}
O.~R. Musin and A.~S. Tarasov.
\newblock The strong thirteen spheres problem.
\newblock {\em Discrete Comput. Geom.}, 48(1):128--141, 2012.

\bibitem{MTT13}
O. R. Musin, A.~S. Tarasov. Enumeration of irreducible contact graphs on the sphere. {\it J. Math. Sci.,} 203:837--850, 2014.

\bibitem{MTT14}
O. R. Musin, A.~S. Tarasov, The Tammes problem for N=14.  {\it Experimental Math.,} 24:460--468, 2015.  

\bibitem{MT15a}
O. R. Musin, A. S. Tarasov, Extreme problems of circle packings on a sphere and irreducible contact graphs. {\it Proc.  Steklov Inst. of Math,}  288:117--131, 2015. 

\bibitem{MusVo}
O. R. Musin and A.\,Yu. Volovikov, Borsuk--Ulam type spaces.  {\it Moscow Math. J.,} {15(4)}:749--766, 2015. 


\bibitem{odl79}
A.~M. Odlyzko and N.~J.~A. Sloane.
\newblock New bounds on the number of unit spheres that can touch a unit sphere
  in {$n$} dimensions.
\newblock {\em J. Combin. Theory Ser. A}, 26(2):210--214, 1979.

\bibitem{pfe04}
F.~Pfender and G.~M. Ziegler.
\newblock Kissing numbers, sphere packings, and some unexpected proofs.
\newblock {\em Notices Amer. Math. Soc.}, 51(8):873--883, 2004.

\bibitem{rob61}
R.~M. Robinson.
\newblock Arrangements of 24 points on a sphere.
\newblock {\em Math. Ann.}, 144:17--48, 1961.

\bibitem{Scho}
I.~J. Schoenberg,  Positive definite functions on spheres, {\em Duke Math. J.},
{9}:96-107, 1942.

%\bibitem{sch38}
%I.~J. Schoenberg.
%\newblock Metric spaces and completely monotone functions.
%\newblock {\em Ann. of Math. (2)}, 39(4):811--841, 1938.

\bibitem{sch05}
A.~Schrijver.
\newblock New code upper bounds from the {T}erwilliger algebra and semidefinite
  programming.
\newblock {\em IEEE Trans. Inform. Theory}, 51(8):2859--2866, 2005.

\bibitem{sch51}
K.~Sch{\"u}tte and B.~L. van~der Waerden.
\newblock Auf welcher {K}ugel haben {$5$}, {$6$}, {$7$}, {$8$} oder {$9$}
  {P}unkte mit {M}indestabstand {E}ins {P}latz?
\newblock {\em Math. Ann.}, 123:96--124, 1951.

\bibitem{sch53}
K.~Sch{\"u}tte and B.~L. van~der Waerden.
\newblock Das {P}roblem der dreizehn {K}ugeln.
\newblock {\em Math. Ann.}, 125:325--334, 1953.

\bibitem{schur09}
A. Sch{\"u}rmann, \emph{Computational geometry of positive definite
  quadratic forms}, University Lecture Series 49, 2009.


%\bibitem{szp03}
%G.~G. Szpiro.
%\newblock {\em Kepler's conjecture}.
%\newblock John Wiley \& Sons Inc., Hoboken, NJ, 2003.
%\newblock How some of the greatest minds in history helped solve one of the
 % oldest math problems in the world.

%\bibitem{tam30}
%P.~M.~L. Tammes.
%\newblock {\em On the origin of number and arrangement of the places of exit on
%  the surface of pollen-grains}.
%\newblock Amsterdam, de Bussy, 1930.

\bibitem{viaz17}
M. Viazovska. The sphere packing problem in dimension 8, {\em Ann. of Math.}, 185:991--1015, 2017. 

\bibitem{zong}
C. Zong, Sphere Packings, Springer-Verlag, New York, 1999.

}
\end{thebibliography}
\end{document}